\documentclass[12pt]{amsart}
\usepackage{amscd,amssymb}
\usepackage[graph,frame,poly,arc]{xy}  
\usepackage[plainpages,backref,urlcolor=blue]{hyperref}
\usepackage{mathtools}

\theoremstyle{plain}
\newtheorem{thm}[subsection]{Theorem}

\newtheorem{prop}[subsection]{Proposition}
\newtheorem{cor}[subsection]{Corollary}

\theoremstyle{definition}
\newtheorem{rk}[subsection]{Remark}

\numberwithin{equation}{section}
\newcommand{\C}{\mathbb{C}}
\newcommand{\PP}{\mathbb{P}}

\DeclarePairedDelimiter\ceil{\lceil}{\rceil}

\begin{document}

\title [On the freeness of rational cuspidal plane curves]
{On the freeness of rational cuspidal plane curves}

\author[Alexandru Dimca]{Alexandru Dimca}
\address{Universit\'e C\^ ote d'Azur, CNRS, LJAD and INRIA, France }
\email{dimca@unice.fr}

\author[Gabriel Sticlaru]{Gabriel Sticlaru}
\address{Faculty of Mathematics and Informatics,
Ovidius University
Bd. Mamaia 124, 900527 Constanta,
Romania}
\email{gabrielsticlaru@yahoo.com }



\begin{abstract} We bring additional support to the conjecture saying that a rational cuspidal plane curve is either free or nearly free. This conjecture was confirmed for curves of even degree, and in this note we prove it for many odd degrees.
In particular, we show that this conjecture holds for  the curves of degree at most 34.

\end{abstract}
 
\maketitle


\section{Introduction} 

A plane rational cuspidal curve is a rational curve $C:f=0$ in the complex projective plane $\PP^2$, having only unibranch singularities. The study of these curves has a long and fascinating history, some long standing conjectures, as the Coolidge-Nagata conjecture being
proved only recently, see \cite{KP0}, other conjectures, as the one on the number of singularities of such a curve being bounded by 4, see \cite{Pion}, are still open.
The classification of such curves is not easy, there are a wealth of examples even when additional strong restrictions are imposed, see \cite{FLMN, F, FZ, Moe,PP, SaTo}.

Free divisors, defined by a homological property of their Jacobian ideals, have been introduced in a local analytic setting by K. Saito in \cite{KS},
and then extended to projective hypersurfaces, see \cite{BC, ST} and the references there.
We have remarked in \cite{DStFD} that many plane rational cuspidal curves are free. The remaining examples of plane rational cuspidal curves in the available classification lists turned out to satisfy a weaker homological property, which was chosen as the definition of a nearly free curve, see \cite{DStRIMS}. Subsequently, a number of authors have establish interesting properties of this class of curves, see \cite{B+, MV}.

In view of the above remark, we have conjectured in \cite[Conjecture 1.1]{DStRIMS} that any plane rational cuspidal curve $C$ is either free or nearly free. This conjecture was proved in \cite[Theorem 3.1]{DStRIMS} for curves $C$ whose degree $d$ is even, as well as for some cases when $d$ is odd, e.g. when $d=p^k$, for a prime number $p>2$. In this note we take a closer look at the case $d$ odd.

Let $S=\C[x,y,z]$ be the polynomial ring in three variables $x,y,z$ with complex coefficients, $f \in S$ a reduced homogeneous polynomial of degree $d \geq 2$, and let $f_x, f_y$ and $f_z$  be the partial derivatives of $f$ with respect to $x,y$ and $z$ respectively.
Consider the graded $S-$submodule $AR(f) \subset S^{3}$ of {\it all relations} involving the derivatives of $f$, namely
$$\rho=(a,b,c) \in AR(f)_q$$
if and only if  $af_x+bf_y+cf_z=0$ and $a,b,c$ are in $S_q$, the space of homogeneous polynomials of degree $q$. 
The minimal degree of a Jacobian relation for the polynomial $f \in S_d$ is the integer $mdr(f)$
defined to be the smallest integer $m\geq 0$ such that $AR(f)_m \ne 0$.
When $mdr(f)=0$, then $C:f=0$ is a union of  lines passing through one point, and hence $C$ is cuspidal only for $d=1$. 

We assume from now on in this note that 
$ mdr(f)\geq 1.$ It turns out that a rational cuspidal curve $C:f=0$ with $mdr(f)=1$ is nearly free. Indeed,  this follows from \cite[Proposition 4.1]{Drcc}. To see this, note that the implication (1) $\Rightarrow $ (2) there  holds for any $d \geq 2$.

Assume from now on that $d$ is odd, and let 
\begin{equation}
\label{pfactor}
d=p_1^{k_1} \cdot p_2^{k_2} \cdots p_m^{k_m}
\end{equation}
be the prime decomposition of $d$. We assume also that $m\geq 2$, the case $m=1$ of our conjecture being settled in \cite[Corollary 3.2]{DStRIMS}. By changing the order of the $p_j$'s if necessary, we can and do assume that  $p_1^{k_1}>p_j^{k_j}$,
for any $2 \leq j \leq m$. Set $e_1=d/p_1^{k_1}.$
With these assumptions and notations, the main results of this note are the following.

\begin{thm}
\label{thmA}
Let $C:f=0$ be a rational cuspidal curve of degree $d=2d'+1$ an odd number. Then $mdr(f) \leq d'$ and if equality holds, then $C$ is either free or nearly free.

\end{thm} 

\begin{thm}
\label{thmB}
Let $C:f=0$ be a rational cuspidal curve of degree $d=2d'+1$,  an odd number as in \eqref{pfactor}. Then, if
$$mdr(f) \leq r_0:=\frac{d-e_1}{2},$$ then $C$ is either free or nearly free. In particular, the following hold.

\begin{itemize}
		\item[i)]  If $d=3p^k$, with $p$ a prime number, then $C$ is either free or nearly free.
		\item[ii)]  $d=5p^k$, with $p$ a prime number, $p^k>3$, then $C$ is either free or nearly free, unless $mdr(f)=d'-1$.
		
	\end{itemize}

\end{thm}

\begin{rk}
\label{rkA}
Note that, for $d \ne 15$, we have $e_1 \leq d/7$ and hence
$$r_0=\frac{d-e_1}{2} \geq \ceil*{\frac{d(1-\frac{1}{7})}{2}}=\ceil*{\frac{3d}{7}}.$$
Therefore, the only cases not covered by our results correspond  to curves of odd degree $d$, such that $r=mdr(f)$ satisfies
$$\ceil*{\frac{3d}{7}}+1 \leq r_0+1 \leq r \leq d'-1=\frac{d-3}{2}.$$
\end{rk}

\begin{cor}
\label{corA}
A rational cuspidal curve $C:f=0$ of degree $d$ is either free or nearly free, if 
one of the following holds.

\begin{enumerate}

\item $mdr(f) \leq 15$, or

\item $d \leq 90$, unless we are in one of the following situations.

\begin{itemize}
		\item[i)]  $d=35$ and $mdr(f)=16$;
		\item[ii)]  $d=45$ and $mdr(f)=21$;
		\item[iii)]  $d=55$ and $mdr(f)=26$;
		\item[iv)]  $d=63$ and $mdr(f)\in \{29,30\}$;
		\item[v)]  $d=65$ and $mdr(f)=31$.
		\item[vi)]  $d=77$ and $mdr(f)\in \{36,37\}$.
		\item[vii)]  $d=85$ and $mdr(f)=41$.
	\end{itemize}

\end{enumerate}

\end{cor}
In the excluded situations, our results do not allow us to conclude.

The proof of our main results are based on a deep result by  U. Walther, see \cite[Theorem 4.3]{Wa} bringing into the picture the monodromy of the Milnor fiber $F:f=1$ associated to the curve $C:f=0$. A second ingredient is our results on the relations between  the Hodge filtration and pole order filtration on the cohomology group $H^1(F,\C)$, see \cite[Theorem 1.2]{DStFor} and \cite[Proposition 2.2]{DStproj}.

\bigskip

The first author thanks AROMATH team at INRIA Sophia-Antipolis for excellent working conditions, and in particular Laurent Bus\'e for stimulating discussions.

\section{Some facts about free and nearly free curves} 

Here we recall some basic notions on free and nearly free curves.  We denote by $J_f$ the Jacobian ideal of $f$, i.e. the homogeneous ideal of $S$ spanned by the partial derivatives $f_x,f_y,f_z$ and let $M(f)=S/J_f$ be the corresponding graded ring, called the Jacobian (or Milnor) algebra of $f$. Let $I_f$ denote the saturation of the ideal $J_f$ with respect to the maximal ideal ${\bf m}=(x,y,z)$ in $S$ and recall the relation with the 0-degree local cohomology 
$$
N(f):=I_f/J_f=H^0_ {\bf m}(M(f)).
$$
It was shown in \cite[Corollary 4.3]{DPop} that the graded $S$-module  $N(f)$ satisfies a Lefschetz type property with respect to multiplication by generic linear forms. This implies in particular the inequalities
\begin{equation}
\label{Lef}
0 \leq n(f)_0 \leq n(f)_1 \leq ...\leq n(f)_{[T/2]} \geq n(f)_{[T/2]+1} \geq ...\geq n(f)_T \geq 0,
\end{equation}
where $T=3d-6$ and $n(f)_k=\dim N(f)_k$ for any integer $k$. 
If we set 
$\nu(f)=\dim N(f)_{[T/2]},$
then $C:f=0$ is a free curve if $\nu(f)=0$.
We say that  $C:f=0$ is a nearly free curve if $\nu(f) =1$, see \cite{B+, Dmax, Drcc, DStFD, DStRIMS} for more details, equivalent definitions and many examples.

Note that the curve $C:f=0$ is  { free } if and only if the  graded $S$-module $AR(f)$ is free of rank 2, i.e. there is an isomorphism of graded $S$-modules
$$AR(f)=S(-d_1) \oplus S(-d_2)$$
for some positive integers $d_1 \leq d_2$. When $C$ is free, the integers $d_1 \leq d_2$  are called the {exponents} of $C$.  They satisfy the relations 
\begin{equation}
\label{free1}
 d_1+d_2=d-1 \text{ and } \tau(C)=(d-1)^2 - d_1d_2,
\end{equation}
where $\tau(C)$ is the total Tjurina number of $C$,  that is
$\tau(C)=\sum_{i=1}^p \tau(C,{x_i})$, the $x_i$'s being the singular points of $C$, and $\tau(C,{x_i})$ denotes the Tjurina number of the isolated plane curve singularity $(C,x_i)$,
see for instance \cite{DS14, DStFD}. In the case of a nearly free curve, there are also the exponents $d_1 \leq d_2$, and this time they verify
\begin{equation}
\label{free1.5}
 d_1+d_2=d \text{ and } \tau(C)=(d-1)^2 - d_1(d_2-1)-1.
\end{equation}
Both for a free and a nearly free curve $C:f=0$, one has $mdr(f)=d_1$, and hence
$mdr(f) \leq (d-1)/2$ for a free curve $C$, and $mdr(f) \leq d/2$ for a nearly free curve $C$. 
 It follows that
Theorem \ref{thmA} gives a similar inequality for any rational cuspidal curve. Our examples of rational cuspidal curves given in \cite{DStExpo}, which are also free or nearly free, show that all the possible values of $mdr(f)$ do actually occur for any fixed degree $d$.
 It follows that

 If we set $r=mdr(f)$, then the curve $C:f=0$ is free (resp. nearly free) if and only if 
\begin{equation}
\label{free2}
\tau(C)=\tau(d,r):=(d-1)^2-r(d-1-r)
\end{equation}
(resp.
$\tau(C)=\tau(d,r)-1$), see \cite{Dmax}.

\begin{rk}
\label{rkB}
If the equation $f=0$ of the curve $C$ is given explicitly, then one can use a computer algebra software, for instance Singular \cite{Sing}, in order to compute the integer $mdr(f)$.
Such a computer algebra software can of course decide whether the curve $C$ is free or nearly free, see for instance the corresponding code on our website\\
\url{http://math.unice.fr/~dimca/singular.html}\\
However, for large degrees $d$, it is much quicker to determine the integer $mdr(f)$.
\end{rk}

\section{The proofs} 

First we recall the setting used in the proof of \cite[Theorem 3.1]{DStRIMS}.
The key results of U. Walther in \cite[Theorem 4.3]{Wa}  yield the inequality
\begin{equation}
\label{UW}
\dim N(f)_{2d-2-j} \leq \dim H^2(F, \C)_{\lambda},
\end{equation}
for $j=1,2,...,d$, where $F: f(x,y,z)-1=0$ is the Milnor fiber in $\C^3$ associated to the plane curve $C$, and the subscript $\lambda$ indicates the eigenspace of the monodromy action corresponding to the eigenvalue $\lambda= \exp(2\pi i (d+1-j)/d)=\exp(-2\pi i (j-1)/d)$. 

Assume that $C$ is a rational cuspidal curve of degree $d$. Denote  by $U$ the complement $\PP^2 \setminus C$, and note that its topological Euler characteristic is given by
$E(U)=E(\PP^2)-E(C)=1$. Since $F$ is a cyclic $d$-fold covering of the complement $U$, 
it follows that $H^m(F,\C)_1=H^m(U, \C)=0$ for $m=1,2$. We have also
$$\dim H^2(F,\C)_{\lambda} -\dim H^1(F,\C)_{\lambda} +\dim H^0(F,\C)_{\lambda}=E(U)=1,$$
see for instance \cite[Prop. 1.21, Chapter 4]{D1} or \cite[Cor. 5.1 and Remark 5.1]{DHA}.
For any $\lambda \ne 1$, since clearly  $H^0(F,\C)_{\lambda}=0$, we get 
\begin{equation}
\label{UW2}
\dim H^2(F,\C)_{\lambda} =\dim H^1(F,\C)_{\lambda}+1. 
\end{equation}

\subsection{Proof of Theorem \ref{thmA}} 
Suppose now that $d$ is odd, say $d=2d'+1$. In order to prove $\nu(f) \leq 1$, in view of the inequality \eqref{UW}, it is enough to show that $\dim H^2(F,\C)_{\lambda} =1$, for  ${\lambda= \exp(2\pi i d'/(2d'+1))   }$ which corresponds to $j=d'+2$.
The equation \eqref{UW2} tells us that this is equivalent to $\dim H^1(F,\C)_{\lambda} =0$.
Using \cite[Proposition 2.2]{DStproj}, see also \cite[Remark 4.4]{DStMFgen}, it follows that
$$\dim H^1(F,C)_{\lambda}= \dim E^{1,0}_2(f)_k + \dim E^{1,0}_2(f)_{d-k},$$
where $k=j-1=d'+1$. Here $E^{1,0}_2(f)_k$ and  $ E^{1,0}_2(f)_{d-k}$ denote some terms of the second page of spectral sequences used to compute the monodromy action on the Milnor fiber $F$, see \cite{D1, DStFor, DStMFgen, Sa3} for details. Note also that the weaker result in \cite[Theorem 1.2]{DStFor}  is enough for this proof.
By the construction of these spectral sequences, it follows that, for $q\leq d $, one has an identification
$$E^{1,0}_2(f)_q=\{(a,b,c) \in AR(f)_{q-2} \  \ :  \  \ a_x+b_y+c_z=0 \},$$
where $a_x$ is the partial derivative of $a$ with respect to $x$ and so on.
It follows that
$$\dim E^{1,0}_2(f)_k + \dim E^{1,0}_2(f)_{d-k} \leq \dim AR(f)_{d'-1}+\dim AR(f)_{d'-2} .$$
If $mdr(f) \geq d'$, it follows that $AR(f)_{d'-1}= AR(f)_{d'-2}=0$, and hence the curve $C$ is either free or nearly free. But this implies that $mdr(f) \leq d'$, as explained in the previous section.

\subsection{Proof of Theorem \ref{thmB}} 

For the reader's convenience, we divide this proof into two steps. 

 \begin{prop}
\label{propA}
With the above notation, we have
$ \dim N(f)_j \leq 1$,
for  any integer $j\leq d-3+r_0$ and for any  integer $j \geq 2d-3-r_0$.
\end{prop} 

\proof Let $t_1=(p_1^{k_1}-1)/2$ and note that
$$d-3+r_0=d-3+e_1t_1 <d-3+\frac{d}{2}=\frac{T}{2}.$$
We apply the inequality \eqref{UW} with $j=d+1-e_1t_1=d+1-r_0$. It follows that
\begin{equation}
\label{UW3}
\dim N(f)_{d-3-r_0} \leq \dim H^2(F, \C)_{\lambda},
\end{equation}
with 
$$\lambda= \exp(2\pi i r_0/d)=\exp\left(\frac{2\pi i t_1}{p_1^{k_1}}\right).$$
Since this eigenvalue has order a prime power, it follows from Zariski's Theorem, see \cite[Proposition 2.1]{AD}, 
that $H^1(F,\C)_{\lambda}=0$. Using \eqref{UW2}, we get 
$ \dim N(f)_j \leq 1$ for $j=d-3+r_0$. The claim for $j=2d-3-r_0$ follows from the fact that
the graded module $N(f)$ enjoys a duality property: $\dim N(f)_j=\dim N(f)_{T-j}$, for any integer $j$, see \cite{Se, SW}. The Lefschetz type property of the graded module $N(f)$, see
\eqref{Lef}, completes the proof of this Proposition.

 \begin{prop}
\label{propB}
A rational cuspidal curve $C:f=0$ of degree $d$ as in \eqref{pfactor}, and such that  $r=mdr(f) \leq r_0$,
is either free or nearly free.

\end{prop} 

\proof
We use the formulas (2.2) and (2.3) from \cite{Dmax} and get the following equality
$$ \dim AR(f)_{d-r-2}- \dim N(f)_{2d-r-3}+\dim AR(f)_{r-2}=$$
$$=3{d-r \choose 2} - {2d-r-1 \choose 2}+\tau(C).$$
For any curve $C:f=0$, it is known that, if $\rho_1 \in AR(f)_r$, then the first relation $\rho \in AR(f)_m$ which is not a multiple of $\rho_1$ occurs in a degree $m \geq d-r-1$, see \cite[Lemma 1.1]{ST}.
It follows that 
$$\dim AR(f)_{d-r-2} = \dim S_{d-2r-2}\cdot \rho_1={d-2r \choose 2}=\frac{(d-2r)(d-2r-1)}{2},$$
for any $r$ such that  $2r \leq d$.
Using the obvious fact that $AR(f)_{r-2}=0$, a direct computation shows that
$$\tau(C)=\tau(d,r)-\dim N(f)_{2d-r-3}.$$
Since $r \leq r_0$, it follows that $2d-r-3 \geq 2d-r_0-3$, and hence $\dim N(f)_{2d-r-3}\leq 1$ by Proposition \ref{propA}. The claim follows now using the characterization of free (resp. nearly free) curves given above in \eqref{free2}.
\endproof

It remains to prove the last claim in Theorem \ref{thmB}.
If $d=3p^k$, we can assume $p^k=2p'+1>3$ and then 
$$ r_0=\frac{d-e_1}{2}=3p'.$$
On the other hand, $d=3p^k=6p'+3=2(3p'+1)+1$, and hence $d'=3p'+1$.
Since $mdr(f) \leq d'$ by Theorem \ref{thmA}, we get either $mdr(f)=d'=3p'+1$, and then we conclude by  Theorem \ref{thmA}, or $mdr(f)\leq d'-1=3p'=r_0$, and then we conclude using Theorem \ref{thmB}.

\subsection{Proof of Corollary \ref{corA}} 
To prove the first claim, we have to consider the minimal possible value of
$$ r_0=\frac{d-e_1}{2}=e_1 \frac{p_1^{k_1}-1}{2},$$
when $d$ is odd, but neither a prime power, nor of the form $3p$, with $p>3$ prime.
First, if $p_1=3$, then $k_1 \geq 2$ and $e_1 \geq 5$, hence $r_0 \geq 20$.
Otherwise, $p_1 \geq 5$, $e_1 \geq 3$, but both cannot be equalites.
 It follows that the minimal values are obtained for $p_1=5^2$ and $e_1 = 3$, or $p_1=7$ and $e_1=5$. In the first case we get $r_0 = 36$, in the second we get $r_0=15$.
To prove the second claim, just use Remark \ref{rkA}.


\begin{thebibliography}{00}


\bibitem{AD}  E. Artal Bartolo, A. Dimca, On fundamental groups of plane curve complements, 
Ann. Univ. Ferrara 61(2015), 255-262.



\bibitem{B+} E. Artal Bartolo, L. Gorrochategui, I. Luengo, A. Melle-Hern\' andez,
On some conjectures about free and nearly free divisors, in: {\it Singularities and Computer Algebra, Festschrift for Gert-Martin Greuel on the Occasion of his 70th Birthday}, pp. 1--19, Springer (2017)

\bibitem{BC} R.-O. Buchweitz, A. Conca,  New free divisors from old, J. Commut. Algebra 5 (2013), no. 1, 17--47.




\bibitem{Sing} { W. Decker, G.-M. Greuel, G. Pfister \and H. Sch{\"o}nemann.} \newblock {\sc Singular} {4-0-1} --- {A} computer algebra system for polynomial computations, available at {http://www.singular.uni-kl.de} (2014).


\bibitem{D1} A.~Dimca, {\em Singularities and topology of hypersurfaces}, Universitext, Springer-Verlag, New York, 1992. 



\bibitem{DHA}  A. Dimca,   {\em Hyperplane Arrangements: An Introduction}, Universitext, Springer, 2017


\bibitem{Dmax}  A. Dimca, Freeness versus maximal global Tjurina number for plane curves, 
Math. Proc. Cambridge Phil. Soc.  163 (2017), 161--172.


\bibitem{Drcc} A. Dimca, On rational cuspidal plane curves, and the local cohomology of Jacobian rings, arXiv:1707.05258.



\bibitem{DPop} A. Dimca, D. Popescu, 
Hilbert series and Lefschetz properties of dimension one almost complete intersections, Comm. Algebra 44 (2016), 4467--4482.




\bibitem{DS14} A. Dimca, E. Sernesi,  Syzygies and logarithmic vector fields along plane curves, Journal de l'\'Ecole polytechnique-Math\'ematiques 1(2014), 247-267.



\bibitem{DStExpo} A. Dimca, G. Sticlaru, On the exponents of free and nearly free projective plane curves, Rev. Mat. Complut. 30(2017), 259--268.

\bibitem{DStFor} A. Dimca, G. Sticlaru, A computational approach to Milnor fiber cohomology, Forum Math. 29 (2017),  831-- 846.

\bibitem{DStFD} A. Dimca, G. Sticlaru, Free divisors and rational cuspidal plane curves, Math. Res. Lett. 24(2017), 1023--1042.



\bibitem{DStRIMS} A. Dimca, G. Sticlaru, Free and nearly free curves vs. rational cuspidal plane curves, Publ. RIMS Kyoto Univ., 54:163--179 (2018).

\bibitem{DStMFgen} A. Dimca, G. Sticlaru, Computing the monodromy and pole order filtration on Milnor fiber cohomology of plane curves, arXiv: 1609.06818.

\bibitem{DStproj}  A.Dimca, G. Sticlaru, Computing Milnor fiber monodromy for projective hypersurfaces, arXiv:1703.07146.




\bibitem{FLMN} J. Fernandez de Bobadilla, I. Luengo, A. Melle-Hernandez, A. Nemethi, Classification of rational unicuspidal projective curves whose singularities have one Puiseux pair,
in: {\it Real and Complex Singularities Sao Carlos 2004}, Trends in Mathematics, Birkhauser 2006, pp. 31-45.



\bibitem{F} T. Fenske, Rational 1- and 2-cuspidal plane curves,  Beitr\" age Algebra
Geom., 40(1999), 309-329.

\bibitem{FZ} H. Flenner, M. Zaidenberg, On a class of rational plane curves, Manuscripta Math. 89(1996), 439 - 459.







\bibitem{KP0} M. Koras, K. Palka, The Coolidge--Nagata conjecture. Duke Math. J. 166 (2017), no. 16, 3085--3145.



 
\bibitem{MV} S. Marchesi, J. Vall\`es, Nearly free curves and arrangements: a vector bundle point of view, arXiv:1712.04867.


\bibitem{Moe} T. K.  Moe, Rational Cuspidal Curves, arXiv:1511.02691 (139 pages, Master thesis, University of Oslo,  2008)



\bibitem{PP} K. Palka, T. Pe\l ka,
 Classification of planar rational cuspidal curves. I. C**-fibrations, 
Proc. London Math. Soc. 115 (2017), 638--692.

\bibitem{Pion} J. Piontkowski, On the number of cusps of rational cuspidal plane
curves, Experiment. Math., 16(2007), 251-255.


\bibitem{KS} K. Saito, Theory of logarithmic differential forms and logarithmic vector fields, J. Fac. Sci. Univ. Tokyo Sect. IA Math. 27 (1980), no. 2, 265-291.

\bibitem{Sa3} M. Saito, Bernstein-Sato polynomials  for projective hypersurfaces with weighted homogeneous isolated singularities, arXiv:1609.04801.

\bibitem{SaTo} F. Sakai, K. Tono, Rational cuspidal curves of type (d,d-2) with one or two cusps, Osaka J. Math. 37(2000), 405-415.




\bibitem{Se} E. Sernesi,  The local cohomology of the jacobian ring, {Documenta Mathematica},  19 (2014), 541-565. 



\bibitem{ST} A. Simis, S.O. Toh\u aneanu, Homology of homogeneous divisors, Israel J. Math. 200 (2014), 449-487.



\bibitem{SW} D. van Straten, T. Warmt,  Gorenstein duality for one-dimensional almost complete intersections--with an application to non-isolated real singularities, Math. Proc.Cambridge Phil. Soc.158(2015), 249--268.

\bibitem{Wa} U. Walther, 
The Jacobian module, the Milnor fiber, and the $D$-module generated by $f^s$,  Invent. Math. 207 (2017),1239--1287.




\end{thebibliography}
\end{document}